\documentclass[12pt,reqno]{amsart}

\usepackage{amsmath,amssymb,amsthm,hyperref,graphicx}
\usepackage{latexsym,mathtools, framed}
\usepackage{color}
\usepackage{enumerate}
\usepackage[margin=2.5cm]{geometry}

\newtheorem{thm}{Theorem}[section]
\newtheorem{lemma}[thm]{Lemma}

\newtheorem{corollary}[thm]{Corollary}
\newtheorem{conj}[thm]{Conjecture}
\newtheorem{question}[thm]{Question}
\newtheorem{setup}[thm]{Set-up}
\theoremstyle{definition}
\newtheorem{defn}[thm]{Definition}

\newtheorem{remark}[thm]{Remark}

\DeclareMathOperator{\Ass}{Ass}
\DeclareMathOperator{\maxAss}{maxAss}

\DeclareMathOperator{\conv}{conv}

\newcommand{\mi}{\mathfrak{m}}
\newcommand{\p}{\mathfrak{p}}
\newcommand{\q}{\mathfrak{q}}
\newcommand{\RR}{\mathbb{R}}
\newcommand{\x}{\mathbf{x}}
\newcommand{\ba}{\mathbf{a}}
\newcommand{\bb}{\mathbf{b}}
\newcommand{\bc}{\mathbf{c}}

\newcommand{\br}{\mathbf{r}}
\newcommand{\C}{\chi}
\newcommand{\Om}{\Omega}
\newcommand{\Lm}{\Lambda}
\newcommand{\supp}[1]{\operatorname{supp}\left(#1\right)}
\newcommand{\lattice}[1]{\mathcal L(#1)}
\newcommand{\ww}{\widehat{\alpha}}
\newcommand{\symbpoly}{\mathcal{Q}}
\newcommand{\N}{\mathbb{N}}

\numberwithin{equation}{section}

\title[On the Waldschmidt constant of monomial ideals]{On the Waldschmidt constant of monomial ideals}

\author{}

\date{}

\begin{document}

\author[Nimu Basak]{Nimu Basak}
\address{Department of Mathematics, Indian Institute of Technology Jammu, Jammu \& Kashmir,     INDIA - 181221.}
\email{2024rma1035@iitjammu.ac.in}		
	
\author[Bijender]{Bijender}
\address{Department of Mathematics, Delhi Technological University, Delhi, INDIA - 110042.}
\email{bijenderdharan445@gmail.com}	

\author[Ajay Kumar]{Ajay Kumar}
\address{Department of Mathematics, Indian Institute of Technology Jammu, Jammu \& Kashmir, INDIA - 181221.}
\email{ajay.kumar@iitjammu.ac.in}

\thanks{}

  	

\begin{abstract}
    We investigate a conjecture proposed by Cooper et al. in 2017 regarding a lower bound for the Waldschmidt constant of monomial ideals. We first show that this conjecture does not hold in general by explicitly computing the Waldschmidt constant for a family of monomial ideals. Also, we identify a class of monomial ideals for which this conjecture holds by proving a containment between symbolic and ordinary powers. For a monomial ideal, we provide a new lower bound for the Waldschmidt constant, which modifies the conjectured lower bound under a certain condition.
\end{abstract}

\maketitle

\section{Introduction}\label{S-Introduction}

The comparison between symbolic and ordinary powers of homogeneous ideals is central to addressing major problems in commutative algebra and algebraic geometry. Because explicit algebraic computation of symbolic powers is exceptionally challenging, research largely focuses on the containment problem, that is, determining when the symbolic power of a homogeneous ideal is contained in the ordinary power. Given a homogeneous ideal $I$ in the polynomial ring $R = K[\mathbb{P}^{n - 1}] = K[x_1, \ldots, x_n]$ and a positive integer $m$, the $m^{\text{th}}$ \textit{symbolic power} of $I$ is defined as
\[I^{(m)} = \bigcap\limits_{\p \in \Ass(I)} \left( I^m R_{\p} \cap R \right),\] where the intersection is taken over all associated primes and $R_{\p}$ is the localization of $R$ at the prime ideal $\p.$ The \textit{containment problem} is to determine all pairs of integers $m$ and $t$ for which the inclusion $I^{(m)} \subset I^t$ holds. In contrast, the reverse containment $I^t \subset I^{(m)}$  holds if and only if $m \le t.$ Ein-Lazarsfeld-Smith \cite{ELS}, and Hochster-Huneke \cite{HH} established that for any homogeneous ideal $I \subset R = K[\mathbb{P}^{n - 1}]$ with $n \ge 2$, the containment $I^{(t(n - 1))} \subset I^t$ holds for all $t \ge 1.$ Further progress in this direction was made by Harbourne by showing that  $I^{(t(n - 1) - n + 2)} \subset I^t$ for all $t \ge 1.$ In addition to this, Harbourne and Huneke \cite{HarbourneHuneke} formulated several conjectures regarding the containment of symbolic powers in ordinary powers. Using the symbolic polyhedron, Cooper et al. \cite{CEHH} obtained several containment results regarding symbolic and ordinary powers of monomial ideals. These containment relations yield lower bounds for the Waldschmidt constant of homogeneous ideals. Recall that the \textit{Waldschmidt constant} of a homogeneous ideal $I \subset R$ is given by
\[\widehat{\alpha}(I) = \lim_{m \to \infty} \frac{\alpha(I^{(m)})}{m},\]
where $\alpha(I) = \min\{d : I_d \ne 0\}.$ Moreover, it is known that $\widehat{\alpha}(I) = \inf\left\{\dfrac{\alpha(I^{(m)})}{m} : m \in \mathbb{N} \setminus \{0\} \right\}$; see \cite[Lemma 2.3.1]{BocciBrian}. The study of this invariant was initiated by Waldschmidt and Skoda \cite{Skoda, Waldschmidt} in connection with problems in complex analysis and Diophantine approximation. Later, Chudnovsky \cite{Chudnovsky1} proposed the following conjecture. Let $I$ be the homogeneous ideal defining a finite set of points in $\mathbb{P}^{n - 1}$ over $\mathbb{C}$, where $n \ge 2.$ Then $\ww(I) \ge \dfrac{\alpha(I) + n - 2}{n - 1}.$ Demailly \cite{Demailly} proposed a generalization of Chudnovsky’s Conjecture, which can be stated as follows: Let $I$ be the defining ideal of a finite set of points in $\mathbb{P}^{n - 1}.$ Then for any $m \ge 1$, one has $\ww(I) \ge \dfrac{\alpha(I^{(m)}) + n - 2}{m + n - 2}.$ Recently, H\`a and Sivakumar \cite{Hasiva} established the validity of this conjecture. For further results concerning symbolic powers and their associated asymptotic invariants, we refer to \cite{BTSHB,BGEH,Bisui2,BN,Bocci, EH79,FLPX, Malara,Nagata2, Zar49}. In the setting of monomial ideals, Cooper et al. \cite{CEHH} proposed an analogous conjecture. 

\begin{conj}\label{C1}
    Let $I \subset R = K[x_1, \ldots, x_n]$ be a monomial ideal of big-height $e.$ Then 
    \[\ww(I) \ge \frac{\alpha( I) + e - 1}{e}.\]
\end{conj} 

In order to study Conjecture \ref{C1} for the class of squarefree monomial ideals, the authors in \cite{CEHH} proved the following containment problem: Let $I \subset R = K[x_1, \ldots, x_n]$ be a squarefree monomial ideal of big-height $e$, and let $\mi = \langle x_1, \ldots, x_n \rangle$ be the unique maximal homogeneous ideal of $R.$ Then $I^{(t(m + e - 1) - e + r)} \subset \mi^{(t - 1)(e - 1) + r - 1}(I^{(m)})^t $ for all positive integers $m, t$ and $r.$ This motivates us to study Conjecture \ref{C1} for the classes of non-squarefree monomial ideals. In Section \ref{S-CounterExample}, we construct a family of monomial ideals for which Conjecture \ref{C1} does not hold.

\begin{thm}[Theorem \ref{Th-CounterExample}, Corollary \ref{C-CounterExample}]
    Let $n = 2s$ be an even positive integer. For $1 \le i \le s$, let $w_i \ge 2$ is an integer, and let $I = \q_1 \cap \cdots \cap \q_s$, where
    \[\q_i = \langle x_1, \ldots, x_{i - 1}, x_i^{w_i}, x_{i + 1}, \ldots, x_s, x_{s + i}\rangle.\]
    If $s \ge 2$ and $w = \sum_{i = 1}^{s} \dfrac{w_i}{w_i - 1}$, then $\ww(I) = \dfrac{w}{w - 1}.$ Moreover, if $\sum_{i = 1}^{s} w_i < \dfrac{1}{2} s(s + 2)$, then $s \ge 3$ and 
    \[\ww(I) < \dfrac{\alpha(I) + e - 1}{e}.\]
\end{thm} 

\vspace{2mm}

In Section \ref{S-ContainmentProblem}, we prove the following containment for a class of monomial ideals. 

\begin{thm}[Theorem \ref{Th-MainTheorem}, Corollary \ref{C-LowerBound}]
    Let $I \subset R$ be a monomial ideal with the irredundant irreducible decomposition $I = \q_1 \cap \cdots \cap \q_s.$ For each $1 \le i \le s$, let 
    \[\q_i = \langle x_j^{w_{i, j}} : j \in E_i \rangle,\] 
    where $E_i \subset [n]$, and for $j \in E_i$, $w_{i, j}$ is a positive integer. Assume that $|E_k \cap E_l| \le 1$ for all $k \ne l.$ Then, for every integer $t \ge 1$, we have the containment $I^{(te)} \subset \mi^{t(e - 1)} I^t.$ Moreover, we have 
    \[\ww(I) \ge \dfrac{\alpha(I) + e - 1}{e}.\]
\end{thm} 
    
Let $f \in R$ be a monomial. The \textit{support} of $f$ is given by $\supp{f} = \{j \in [n] : x_j \text{ divides } f\}.$ If $1 \le j \le n$, then the exponent of $x_j$ in $f$ is denoted by $\deg_{x_j}(f).$ Let $\mathcal{G}(I)$ denotes the unique minimal set of monomial generators of a monomial ideal $I \subset R$, and let 
\[w_0 = \min\{\deg_{x_j}(f) : f \in \mathcal{G}(I) \text{ and } j \in \supp f\}.\] 
We next establish the following containment for a monomial ideal $I \subset R.$

\begin{thm}[Theorem \ref{Th-RevisedBound}, Corollary \ref{C-RevisedBound}]
    Let $I \subset R$ be a monomial ideal of big-height $e.$ Then the containment 
    $I^{(t(e + 1))} \subset \mi^{[tw_0]} I^{[t]}$ holds for all $t \ge 1.$ In particular, if $w_0 \ge e - 1$, then 
    \[\ww(I) \ge \dfrac{\alpha(I) + e - 1}{e + 1}.\]
\end{thm}

\section{Preliminaries}

In this section, we recall some basic definitions and notations. We continue to employ the definitions and notations given in Section \ref{S-Introduction}.

Let $\N$ denote the set of non-negative integers, and let $\RR_{+}$ denote the set of non-negative real numbers. Let $\ba = (a_1, \ldots, a_n) \in \N^n.$ We define $|\ba| = \sum_{j = 1}^{n} a_j.$ If $\bb = (b_1, \ldots, b_n) \in \N^n$ satisfy $a_j \le b_j$ for all $1 \le j \le n$, then we write $\ba \le \bb.$

Let $K$ be a field. We write $R$ to denote the polynomial ring $K[x_1,\ldots,x_n].$ The unique homogeneous maximal ideal $\langle x_1, \ldots, x_n \rangle$ of $R$ is denoted by $\mi.$ If $\ba = (a_1, \ldots, a_n) \in \N^n$, then we write $\x^{\ba}$ to denote the monomial $x_1^{a_1} \cdots x_n^{a_n} \in R.$ The \textit{degree} of a monomial $f = \x^{\ba} \in R$ is given by $\deg(f) = |\ba|.$

\begin{defn}
    Let $I \subset R$ be a monomial ideal. 
    \begin{enumerate}[(a)]
        \item A prime ideal $\p \subset R$ is called an \textit{associated prime} ideal of $I$ if there exists $f \in R$ such that $I : f = \p.$ The set of all associated prime ideals of $I$ is denoted by $\Ass(I)$, and the set of all maximal elements of $\Ass(I)$ under inclusion is denoted by $\maxAss(I).$

        \item The \textit{big-height} of $I$ is the maximum of the heights of its associated primes.

    \end{enumerate}
\end{defn}

\begin{defn}
    Let $I \subset R$ be a monomial ideal. 
    \begin{enumerate}[(a)]
        \item An ideal $\q \subset R$ is said to be \textit{irreducible} if $\q = \q_1 \cap \q_2$ for monomial ideals $\q_1, \q_2 \subset R$, then either $\q = \q_1$ or $\q = \q_2.$

        \item A presentation of $I$ as the intersection $I = \q_1 \cap \cdots \cap \q_s$ of irreducible monomial ideals $\q_1, \ldots, \q_s$ is called an \textit{irreducible  decomposition} of $I.$ An irreducible decomposition $I = \q_1 \cap \cdots \cap \q_s$ is called \textit{irredundant} if $\bigcap_{i \ne k} \q_i \not\subset \q_k$ for all $1 \le k \le s.$ 
    \end{enumerate}
\end{defn}

\begin{remark}\label{R-IrreducibleIdeals}
    By \cite[corollary 1.3.2]{HBook}, a monomial ideal $\q \subset R$ is irreducible if and only if there exists $E \subset [n]$ such that $\q = \langle x_{j}^{w_j} : j \in E \rangle$, where for each $j \in E$, $w_j \ge 1$ is an integer.
\end{remark}

We conclude this section by introducing the notion of Frobenius powers.

\begin{defn}
     Let $I \subset R$ be a monomial ideal, let $\mathcal{G}(I)$ be the unique minimal set of monomial generators of $I$, and let $t$ be a positive integer. Then the $t$-th \textit{Frobenius power} $I^{[t]}$ of $I$ is defined as 
     \[I^{[t]} = \langle f^t : f \in \mathcal{G}(I) \rangle.\]
\end{defn}

\section{On a Counterexample to Conjecture \ref{C1}}\label{S-CounterExample}

The main aim of this section is to introduce a class of monomial ideals which disproves Conjecture \ref{C1}. We explicitly compute the Waldschmidt constant for each member of this class by using the notion of symbolic polyhedron, introduced in \cite{CEHH}.
    
\begin{defn}
    Let $A \subset \RR^n$, and let $I\subset R$ be a monomial ideal with the irredundant irreducible decomposition $I = \q_1 \cap \cdots \cap \q_s.$
    \begin{enumerate}[(a)]
        \item A \textit{convex combination} of points $\ba_1, \ldots, \ba_k \in \RR^n$ is the sum $\sum_{i = 1}^{k} \lambda_i \ba_i$ with $\lambda_i \in \RR_{+}$ and $\sum_{i = 1}^{k} \lambda_i = 1.$
        
        \item We say that $A$ is a \textit{convex set} if for any two points $\ba, \bb \in A$, we have $\lambda \ba + (1 - \lambda) \bb \in A$ for all $\lambda \in \RR_{+}$ with $0 \le \lambda \le 1.$

        \item The \textit{convex hull} of $A$, denoted by $\conv(A)$, is the convex set of all finite convex combinations of points of $A.$ Equivalently, it is the smallest convex set containing $A.$
    
        \item The \textit{set of exponent vectors} of the monomials in $I$ is given by
        \[\lattice{I} = \{\ba \in \mathbb{N}^n : \x^{\ba} \in I\}.\]
        Similarly, if $\mathcal{G}(I)$ is the unique minimal set of monomial generators of $I$, then $\lattice{\mathcal{G}(I)}$ denotes the set of exponent vectors of the minimal monomial generators of $I.$

        \item For each $\p \in \maxAss(I)$, we set $\q_{\subset \p} = R \cap IR_{\p}$, where $R_{\p}$ is the localization of $R$ at $\p.$ Following \cite[Definition~5.3]{CEHH}, we define the \textit{symbolic polyhedron} of $I$ by
        \[\symbpoly(I) = \bigcap_{\p \in \maxAss(I)} \conv \left( \lattice{\q_{\subset \p}} \right).\]

     \end{enumerate}
\end{defn}

The following lemma characterize the symbolic polyhedron for irreducible monomial ideals.

\begin{lemma}\label{L-SymbolicPolyhedron}
    Let $E \subset [n]$, and let $I = \langle x_{j}^{w_j} : j \in E \rangle \subset R$ be an irreducible monomial ideal, where for each $j \in E$, $w_j \ge 1$ is an integer. Then 
    \[\conv(\lattice{I}) = \left\{ \ba = (a_1, \ldots, a_n) \in \RR_{+}^{n} : \sum_{j \in E} \dfrac{a_j}{w_j} \ge 1 \right\}.\]
\end{lemma} 

\begin{proof}
    For simplicity, write 
    \[W = \left\{ \ba = (a_1, \ldots, a_n) \in \RR_{+}^{n} : \sum_{j \in E} \dfrac{a_j}{w_j} \ge 1 \right\}.\]
   
    Let $\ba = (a_1, \ldots, a_n) \in \conv(\lattice{I}).$ Then, by applying \cite[Lemma 5.2]{CEHH}, there exist
    $\bb = (b_1, \ldots, b_n) \in \conv(\lattice{\mathcal{G}(I)})$ and $\bc = (c_1, \ldots, c_n) \in \RR_{+}^n$ such that $\ba = \bb + \bc.$ Note that if $\mathbf{e}_j$ denotes the $j$-th standard basis vector of $\RR^{n}$, then 
    \[\lattice{\mathcal{G}(I)} = \{w_j\mathbf{e}_j : j \in E\}.\] 
    This implies that $\bb = \sum_{j \in E} \lambda_j w_j \mathbf{e}_j$, where $\lambda_j \in \RR_{+}$ for each $j \in E$ and $\sum_{j \in E} \lambda_j = 1.$ By comparing the components of vector $\ba$, we obtain $a_j = \lambda_j w_j + c_j$ if $j \in E$, and $a_j = c_j$ if $j \notin E.$ Since $\bc \in \RR_{+}^n$, $\lambda_j \in \RR_{+}$ and $w_j \ge 1$, it follows that $\ba \in \RR_{+}^{n}$ and
    \[\sum_{j \in E} \dfrac{a_j}{w_j} = \sum_{j \in E} \left( \lambda_j + \dfrac{c_j}{w_j} \right) \ge \sum_{j \in E} \lambda_j = 1.\]
    This implies $\ba \in W.$  

    To prove the reverse inclusion, let $\ba = (a_1, \ldots, a_n) \in W$, and let $\lambda = \sum_{j \in E} \dfrac{a_j}{w_j} \ge 1.$ For $j \in E$, let $\lambda_j = \dfrac{a_j}{\lambda w_j}.$ Then $\lambda_j \ge 0$ and $\sum_{j \in E} \lambda_j = 1$, and hence $\sum_{j \in E} \lambda_j w_j \mathbf{e}_j \in \conv(\lattice{\mathcal{G}(I)}).$ In addition to this, since $\ba \in \RR_{+}^{n}$ and $\lambda \ge 1$, we obtain the vector inequality
    \[\ba = \sum_{j = 1}^{n} a_j \mathbf{e}_j \ge \sum_{j \in E} a_j \mathbf{e}_j = \lambda\left( \sum_{j \in E} \lambda_j w_j \mathbf{e}_j \right) \ge \sum_{j \in E} \lambda_j w_j \mathbf{e}_j.\] 
    This implies that $\ba - \sum_{j \in E} \lambda_j w_j \mathbf{e}_j \in \RR_{+}^n.$ Consequently, $\ba \in \conv(\lattice{I}).$
\end{proof}

\begin{remark}\label{R-CramerRule} 
    Consider the $s \times s$ non-homogeneous linear system in variables $X_1, \ldots, X_s$ given by $M X = B$, where
    \[M = 
    \begin{bmatrix}
        1 & w_1 & \cdots & w_1 \\
        w_2 & 1 & \cdots & w_2 \\
        \vdots & \vdots & \ddots & \vdots \\
        w_s & w_s & \cdots & 1
    \end{bmatrix},
    X = 
    \begin{bmatrix}
        X_1 \\
        X_2 \\
        \vdots \\
        X_s
    \end{bmatrix},
    \text{ and } B = 
    \begin{bmatrix}
        w_1 \\
        w_2 \\
        \vdots \\
        w_s
    \end{bmatrix},\] and $w_i \ge 2$ for all $i.$
    Then the determinant of the coefficient matrix $M$ is given by
    \[\det(M) = (1 - w) \prod_{i = 1}^{s} (1 - w_i),\]
    where $w = \sum_{i = 1}^{s}\dfrac{w_i}{w_i - 1}.$ For each $1 \le i \le s$, let $M_i$ denote the matrix obtained by replacing the $i$-th column of the coefficient matrix $M$ with the constant column vector $B.$ Then the determinant of $M_i$ is given by
    \[\det(M_i) = w_i \prod_{\substack{j = 1 \\ j \ne i}}^{s} (1 - w_j).\]
    Thus, by Cramer's rule, the solution of the linear system $MX = B$ is given by
    \[X_i = \dfrac{w_i}{(w_i - 1)(w - 1)}, \qquad \text{ where } 1 \le i \le s.\] 
\end{remark}

The following theorem is the main result of this section.

\begin{thm}\label{Th-CounterExample}
    Let $n = 2s$ be an even positive integer. For $1 \le i \le s$, let $w_i \ge 2$ is an integer, and let $I = \q_1 \cap \cdots \cap \q_s$, where
    \[\q_i = \langle x_1, \ldots, x_{i - 1}, x_i^{w_i}, x_{i + 1}, \ldots, x_s, x_{s + i}\rangle.\]
    If $s \ge 2$ and $w = \sum_{i = 1}^{s} \dfrac{w_i}{w_i - 1}$, then $\ww(I) = \dfrac{w}{w - 1}.$
\end{thm}

\begin{proof}
   In view of \cite[Corollary 6.3]{CEHH}, we have $\ww(I) = \min \{\ba \cdot \mathbf{1} : \ba \in \symbpoly(I)\}$, where $\mathbf{1}$ is the column vector in $\RR_{+}^{n}$ with each entry equal to $1.$ We first show that $\ww(I) \le \dfrac{w}{w - 1}.$ For this, we find a vector $\bb \in \symbpoly(I)$ such that and $\bb \cdot \mathbf{1} = \dfrac{w}{w - 1}.$ Consider the $s \times s$ non-homogeneous linear system in variables $X_1, \ldots, X_s$ given by
    \begin{equation}\label{EQN-LinearSystem}
        \begin{aligned}
            X_1 + w_1 X_2 + \cdots + w_1 X_s & = w_1,\\
            w_2 X_1 + X_2 + \cdots + w_2 X_s & = w_2,\\
            &\ \vdots \\
            w_s X_1 + w_s X_2 + \cdots + X_s & = w_s.
        \end{aligned}
    \end{equation}
    By the definition of the symbolic polyhedron, we have $\symbpoly(I) = \bigcap_{i = 1}^{s} \conv(\lattice{\q_i})$, while it follows from Lemma \ref{L-SymbolicPolyhedron} that
    \[\conv(\lattice{\q_i}) = \Bigg\{ \ba = (a_1, \ldots, a_n) \in \RR_{+}^{n} : \dfrac{a_i}{w_i} + \sum_{\substack{j = 1 \\ j \ne i}}^{s} a_j + a_{s + i} \ge 1 \Bigg\}\]
    for all $1 \le i \le s.$ Thus, if any vector $(a_1, \ldots, a_s) \in \RR_{+}^{s}$ is a solution of the linear system \ref{EQN-LinearSystem}, then the vector $(a_1, \ldots, a_s, 0, \ldots, 0) \in \symbpoly(I).$ By Remark \ref{R-CramerRule}, the solution of the linear system \ref{EQN-LinearSystem} is given by
    \begin{equation}\label{EQN-Solution}
        X_i  = \dfrac{w_i}{(w_i - 1)(w - 1)},  \text{ where } 1 \le i \le s \text{ and } w = \sum_{i = 1}^{s}\dfrac{w_i}{w_i - 1}. 
    \end{equation}
    Consequently, $\bb = (b_1, \ldots, b_s, 0, \ldots, 0) \in \symbpoly(I),$ where $b_i=X_i, 1 \le i \le s.$ It is immediate that $\bb \cdot \mathbf{1} = \sum_{i = 1}^{s} b_i = \dfrac{w}{w - 1}.$ Thus, we obtain $\ww(I) \le \dfrac{w}{w-1}.$ 
    
    To prove the reverse inequality, it is enough to prove that $\ba \cdot \mathbf{1} \ge \dfrac{w}{w - 1}$ for all $\ba \in \symbpoly(I).$ For this, let $\ba = (a_1, \ldots, a_n) \in \symbpoly(I).$ Then, for all $1 \le i \le s$, we have
    \[\dfrac{a_i}{w_i} + \sum_{\substack{j = 1 \\ j \ne i}}^{s} a_j + a_{s + i} \ge 1.\]
    Multiplying by $b_i = \dfrac{w_i}{(w_i - 1)(w - 1)}$ and summing yields
    \begin{equation}\label{INQ-1}
        \sum_{i = 1}^{s} \Bigg( \dfrac{a_i}{w_i} + \sum_{\substack{j = 1 \\ j \ne i}}^{s} a_j \Bigg) b_i + \sum_{i = 1}^{s} a_{s + i} b_i \ge \sum_{i = 1}^{s} b_i = \dfrac{w}{w - 1}.
    \end{equation}
    Now, observe that
    \begin{align*}
        \sum_{i = 1}^{s} \Bigg( \dfrac{a_i}{w_i} + \sum_{\substack{j = 1 \\ j \ne i}}^{s} a_j \Bigg) b_i & = \sum_{i = 1}^{s} \left( \dfrac{a_i b_i}{w_i} - a_i b_i \right) + \sum_{i = 1}^{s} \sum_{j = 1}^{s} a_j b_i \\ 
        & = \sum_{i = 1}^{s} \left( \dfrac{a_i b_i}{w_i} - a_i b_i \right) + \sum_{i = 1}^{s} \sum_{j = 1}^{s} a_i b_j \\ 
        & = \sum_{i = 1}^{s} \dfrac{a_i b_i}{w_i} + \sum_{i = 1}^{s} \sum_{\substack{j = 1 \\ j \ne i}}^{s} a_i b_j \\ 
        & = \sum_{i = 1}^{s} \Bigg( \dfrac{b_i}{w_i} + \sum_{\substack{j = 1 \\ j \ne i}}^{s} b_j \Bigg) a_i. 
    \end{align*}
    Thus, it follows from Inequality \ref{INQ-1} that
    \begin{equation}\label{INQ-2}
        \sum_{i = 1}^{s} \Bigg( \dfrac{b_i}{w_i} + \sum_{\substack{j = 1 \\ j \ne i}}^{s} b_j \Bigg) a_i + \sum_{i = 1}^{s} a_{s + i} b_i \ge \dfrac{w}{w - 1}.
    \end{equation}
   Since $b_i=\dfrac{w_i}{(w_i - 1)(w - 1)}$ and $\sum_{j = 1}^{s} b_j=\dfrac{w}{w - 1},$ we obtain
    \[\dfrac{b_i}{w_i} + \sum_{\substack{j = 1 \\ j \ne i}}^{s} b_j = b_i \left( \dfrac{1}{w_i} - 1 \right) +  \sum_{j = 1}^{s} b_j = \dfrac{b_i (1 - w_i)}{w_i} + \dfrac{w}{w - 1} = 1.\]
    The fact that $w_j \ge 2$ gives $\dfrac{w_j}{w_j - 1} > 1.$ For $1 \le i \le s$, since $s \ge 2$, we have 
    \[w - 1 = - 1 + \sum_{\substack{j = 1 \\ j \ne i}}^{s} \dfrac{w_j}{w_j - 1} + \dfrac{w_i}{w_i - 1} > \dfrac{w_i}{w_i - 1} > 0.\] Therefore, $0<b_i=\dfrac{w_i}{(w_i - 1)(w - 1)}<1.$ Applying these observations to Inequality \ref{INQ-2}, we obtain 
    \[\ba \cdot \mathbf{1} = \sum_{i = 1}^{s} a_i + \sum_{i = 1}^{s} a_{s + i} \ge \sum_{i = 1}^{s} \Bigg( \dfrac{b_i}{w_i} + \sum_{\substack{j = 1 \\ j \ne i}}^{s} b_j \Bigg) a_i + \sum_{i = 1}^{s} a_{s + i} b_i \ge \dfrac{w}{w - 1}.\]
    This completes the proof.
\end{proof}

\begin{remark}
    With the notation of Theorem \ref{Th-CounterExample}, if $s = 2,$ then the fact $\dfrac{1}{w_1 - 1} + \dfrac{1}{w_2 - 1} \le 2$ implies that $w \le 4$. Since $\alpha(I) = 2$ and the big-height $e$ of $I$ is 3,  we obtain
    \[\ww(I) = \dfrac{w}{w - 1} \geq \dfrac{4}{3}=\dfrac{\alpha(I) + e - 1}{e}.\] Therefore, Conjecture \ref{C1} holds in this case.
\end{remark}

As an application of Theorem \ref{Th-CounterExample}, we show that Conjecture \ref{C1} does not hold for $s \geq 3$ for a suitable choice of $w_i$'s. 
\begin{corollary}\label{C-CounterExample}
   With the notation of Theorem \ref{Th-CounterExample}, assume that $\sum_{i = 1}^{s} w_i < \dfrac{1}{2} s(s + 2).$ If $e$ denotes the big-height of $I,$ then 
   \begin{enumerate}[\rm(i)]
       \item $s \ge 3$  and $\ww(I) < \dfrac{\alpha(I) + e - 1}{e}.$
       \item $I^{(te)} \not\subset \mi^{t(e - 1)} I^{t}$  for infinitely many $t \ge 1.$
   \end{enumerate}
\end{corollary}

\begin{proof}
    (i) If $s = 2$, then $w_1 + w_2 \ge 4 = \dfrac{1}{2} s(s + 2)$, since $w_1, w_2 \ge 2.$ This contradicts the hypothesis. Thus $s \ge 3.$ Note that $\alpha(I) = 2$ and $e = s + 1.$ This implies that 
    \[\dfrac{\alpha(I) + e - 1}{e} = \dfrac{s + 2}{s + 1}.\]
    By using Theorem \ref{Th-CounterExample}, we obtain
    \[\ww(I) = \dfrac{w}{w - 1} = 1 + \dfrac{1}{w - 1} \] We claim that $1 + \dfrac{1}{w - 1} < \dfrac{s + 2}{s + 1},$ or equivalently, $w - 1 > s + 1.$
    For simplicity, let $v_i = w_i - 1$, where $1 \le i \le s.$ By using the Cauchy--Schwarz inequality for the vectors $\left( \sqrt{v_1}, \ldots, \sqrt{v_s} \right)$ and $\left( \dfrac{1}{\sqrt{v_1}}, \ldots, \dfrac{1}{\sqrt{v_s}} \right)$, we obtain
    \[\left( v_1 + \cdots + v_s \right) \left( \dfrac{1}{v_1} + \cdots + \dfrac{1}{v_s} \right) \ge s^2.\]
    Since $\sum_{i = 1}^{s} v_i = \sum_{i = 1}^{s} w_i - s < \dfrac{s^2}{2}$, it follows that $\sum_{i = 1}^{s} \dfrac{1}{v_i} = \sum_{i = 1}^{s} \dfrac{1}{w_i - 1} > 2.$ This implies  that
    \[w - 1 = \sum_{i = 1}^{s} \dfrac{w_i}{w_i - 1} - 1 = \sum_{i = 1}^{s} \dfrac{1}{w_i - 1} + s - 1 > s + 1.\]
    (ii) Note that the containment $I^{(te)} \subset \mi^{t(e - 1)} I^{[t]}$ with $t \ge 1$ implies that $\ww(I) \ge \dfrac{\alpha(I) + e - 1}{e}.$  Thus, the result follows from (i). 
\end{proof}

\section{Containment for some special classes}\label{S-ContainmentProblem}

Let $I \subset R$ be a monomial ideal of big-height $e$, and let $I = \q_1 \cap \cdots \cap \q_s$ be the irredundant irreducible decomposition of $I.$ The main aim of this section is to prove the containment 
\[I^{(te)} \subset \mi^{t(e - 1)} I^{[t]}, \qquad t \ge 1,\]
with the condition $|E_k \cap E_l| \le 1$ for all $k \ne l$, where $\q_i = \langle x_j^{w_{i, j}} : j \in E_i \rangle$ for all $1 \le i \le s.$

\begin{setup}\label{setup}
    Let $I \subset R$ be a monomial ideal, and let $I = \q_1 \cap \cdots \cap \q_s$ be the irredundant irreducible decomposition of $I.$ By Remark \ref{R-IrreducibleIdeals}, for each $1 \le i \le s$, we have 
    \[\q_i = \langle x_j^{w_{i, j}} : j \in E_i \rangle,\]
    where $E_i \subset [n]$, and for $j \in E_i$, $w_{i, j}$ is a positive integer. For $1 \le i \le s$, we write $\p_i = \sqrt{\q_i}.$
\end{setup}
    
\begin{lemma}\label{L-SymbolicPowerInequality}
    With the notation of Set-up \ref{setup}, $I^{(m)} \subset \q_1^m \cap \cdots \cap \q_s^m$ for all $m \ge 1.$ Morover, if $\x^{\ba} = x_1^{a_1} \cdots x_n^{a_n} \in I^{(m)}$, then for each $1 \le i \le s$, we have the inequality
    \[\sum_{j \in E_i} \dfrac{a_j}{w_{i,j}} \ge \sum_{j \in E_i} \left\lfloor \dfrac{a_j}{w_{i,j}} \right\rfloor \ge m.\]
\end{lemma}

\begin{proof}
    It follows from \cite[Theorem 3.7]{CEHH} that 
    \[I^{(m)} = \bigcap_{\p \in \maxAss(I)} \q_{\subset \p}^m,\]
    where $\maxAss(I)$ is the set of all maximal associated primes of $I$, and for any $\p \in \maxAss(I)$, we have 
    \[\q_{\subset \p} = \bigcap_{\p_i \subset \p} \q_i.\] 
    Let $1 \le i \le s.$ Then there exists $\p \in \maxAss(I)$ such that $\p_i = \sqrt{\q_i} \subset \p.$ This implies that $\q_{\subset \p} \subset \q_i$, and hence $I^{(m)} \subset \q_{\subset \p}^m \subset \q_i^m.$ Thus, we obtain $I^{(m)} \subset \bigcap_{i = 1}^{s} \q_i^m.$
    
    Now, let $\x^{\ba} \in I^{(m)}$ be a monomial. Then $\x^{\ba} \in \q_i^m.$ Choose a monomial generator $f$ of $\q_i^m$ such that $f$ divides $\x^{\ba}.$ We can write 
    \[f = \prod_{j \in E_i} x_j^{w_{i, j} m_{i, j}}\]
    for some non-negative integers $m_{i, j}$ with $\sum_{j \in E_i} m_{i, j} = m.$ This implies that $w_{i, j} m_{i, j} \le a_j$ for all $j \in E_i$, and hence
    \[\sum_{j \in E_i} \dfrac{a_j}{w_{i,j}} \ge \sum_{j \in E_i} \left\lfloor \dfrac{a_j}{w_{i,j}} \right\rfloor \ge \sum_{j \in E_i} m_{i, j} = m.\]
\end{proof}

\begin{lemma}\label{L-Tool01}
    With the notation of Set-up \ref{setup}, let $\x^{\bb} = x_1^{b_1} \cdots x_n^{b_n} \in I$ and $\x^{\bc} = x_1^{c_1} \cdots x_n^{c_n} \in R$ be such that $\x^{\bc}$ divides $\x^{\bb}.$ Assume that
    \begin{enumerate}[\rm(i)]
        \item $|E_k \cap E_l| \le 1$ for all $k \ne l$; and

        \item there exists $1 \le k \le s$ such that $\supp{\x^{\bc}} \subset E_k$ and $g = \x^{\bb - \bc} \in \q_k.$
    \end{enumerate}
    If $\dfrac{\x^{\bb}}{x_j^{c_j}} \in I$ for all $j \in \supp{\x^{\bc}}$, then $g \in I.$
\end{lemma}

\begin{proof}
    It suffices to show that $g \in \q_i$ for all $i \ne k.$ Let $1 \le i \le s$ with $i \ne k.$ Then, we have
    \[|\supp{\x^\bc} \cap E_i| \le |E_k \cap E_i| \le 1.\] 
    We consider the following two cases.

    \medskip\noindent
    \textit{Case 1.} Assume that $\supp{\x^{\bc}} \cap E_i = \emptyset.$ Since $\x^\bb \in \q_i$, there exists $j \in E_i$ such that $b_j \ge w_{i, j}.$ This implies that $j \notin \supp{\x^\bc}$, and hence $c_j = 0.$ It follows that 
    \[\deg_{x_j}(\x^{\bb - \bc}) = b_j - c_j = b_j \ge w_{i, j}.\]
    Thus, we obtain $g \in \q_i.$ 

    \medskip\noindent
    \textit{Case 2.} Assume that $\supp{\x^{\bc}} \cap E_i = \{p\}.$ Since $\dfrac{\x^{\bb}}{x_p^{c_p}} \in \q_i$, there exists $j \in E_i$ such that $x_j^{w_{i, j}}$ divides $\dfrac{\x^{\bb}}{x_p^{c_p}}.$ If $j = p$, then $b_j - c_j \ge w_{i, j}.$ Now, suppose that $j \ne p.$ The fact $x_j^{w_{i, j}}$ divides $\dfrac{\x^{\bb}}{x_p^{c_p}}$ implies that $b_j \ge w_{i,j}.$ Also, since $\supp{\x^{\bc}} \cap E_i = \{p\}$, we obtain $j \notin \supp{\x^{\bc}}$, and hence $c_j = 0.$ Therefore, $b_j-c_j=b_j \ge w_{i,j.}$ Thus, in both scenarios, we have $b_j - c_j \ge w_{i, j}.$ This implies that $x_j^{w_{i, j}}$ divides $g$, and hence $g \in \q_i.$ 
\end{proof}

\begin{defn}
    With the notation of Set-up \ref{setup}, let $\x^{\bb} = x_1^{b_1} \cdots x_n^{b_n} \in I.$ We define
    \begin{enumerate}[(a)]
        \item $\C(\bb) = \{j \in \supp{\x^{\bb}} : \x^{\bb} / x_j \in I\} = \{j \in \supp{\x^{\bb}} : \x^{\bb} / x_j \in \q_i \text{ for all } i\}$; 

        \vspace{2mm}

        \item $\Om(\bb) = \{j \in \supp{\x^{\bb}} : j \notin \C(\bb)\}$;

        \vspace{2mm}

        \item $\Om(\bb, j) = \{i : 1 \le i \le s \text{ and } \x^{\bb} / x_j \notin \q_i\}$, where $j \in \Om(\bb)$; and

        \vspace{2mm}

        \item $\Lm_i(\bb) = \{j \in E_i : b_j \ge w_{i, j}\}$, where $1 \le i \le s.$
    \end{enumerate}
\end{defn}

\begin{remark}\label{R-MonomialOfI}
    With the notation of Set-up \ref{setup}, let $\x^{\bb} = x_1^{b_1} \cdots x_n^{b_n} \in I.$
    \begin{enumerate}[(i)]
        \item If $\Lm_i(\bb) = \{p\}$ for some $1 \le i \le s$, then $b_p \ge w_{i, p}$ and $b_j < w_{i, j}$ for all $j \in E_i$ with $j \ne p.$ Moreover, if $p \in \C({\bb})$, then $\dfrac{\x^{\bb}}{x_p} \in \q_i.$ This implies that $b_p - 1 \ge w_{i, p}.$
        
        \item For each $1 \le i \le s$, we have $\Lm_i(\bb) \ne \emptyset$, and if $j \in \Lm_i(\bb)$, then $x_j^{b_j} \in \q_i.$ 
        
        \item Let $j \in \Om(\bb)$, and let $i \in \Om(\bb, j).$ Since $\x^{\bb} \in \q_i$, it follows that $\Lm_i(\bb) = \{j\}$ and $b_j = w_{i, j}.$
    \end{enumerate}
\end{remark}

\begin{remark}\label{R-Observation}
    With the notation of Set-up \ref{setup}, let $e$ denote the big-height of $I$, let $t \ge 1$ be an integer, and let $\x^{\ba} = x_1^{a_1} \cdots x_n^{a_n} \in I^{(te)}.$ 
    \begin{enumerate}[(i)]
        \item Let $1 \le i \le s.$ Then, it follows from Lemma \ref{L-SymbolicPowerInequality} that 
        \begin{equation}\label{INQ-3}
            \sum_{j \in E_i} \left\lfloor \dfrac{a_j}{w_{i, j}} \right\rfloor \ge te.
        \end{equation}
        Since $e = \max\{|E_k| : 1 \le k \le s\}$, by the pigeonhole principle, there exists $j_i \in E_i$ such that $\left\lfloor \dfrac{a_{j_i}}{w_{i, j_i}} \right\rfloor \ge t.$ For each $1 \le j \le n$, we write $a_j = t b_j + r_j$, where $0 \le r_j < t.$ Then we have $b_{j_i} + \dfrac{r_{j_i}}{t} \ge w_{i, j_i}.$ Since $b_{j_i}$ and $w_{i, j_i}$ are integers, and $0 \le r_{j_i} < t$, we have $b_{j_i} \ge w_{i, j_i}.$ Consequently, $x_{j_i}^{b_{j_i}} \in \q_i$, and hence $\x^{\bb} \in I$, where $\bb = (b_1, \ldots, b_n).$

        \item Let $f \in \mathcal{G}(I)$ be a monomial generator of smallest degree such that $f$ divides $\x^{\bb}$, and let $\gamma = |\bb| - \deg(f).$ We have 
        \[\x^{\ba} = (\x^{\bb})^t \x^{\mathbf{r}} = f^t \left( \dfrac{\x^{\bb}}{f} \right)^t \x^{\br},\]
        where $\br = (r_1, \ldots, r_n).$ This implies that $f^t$ divides $\x^{\ba}$, and 
        \begin{equation}\label{INQ-4}
            \deg \left( \dfrac{\x^{\ba}}{f^t} \right) = t(|\bb| - \deg(f)) + |\br| = t \gamma + |\br|.  
        \end{equation}
        Since $f^t \in I^{[t]}$, to prove containment $I^{(te)} \subset \mi^{t(e - 1)} I^{[t]}$, we need to establish the inequality $t \gamma + |\br| \ge t(e - 1).$
    \end{enumerate}
\end{remark}    
    
We now present a sequence of lemmas to prove Theorem \ref{Th-MainTheorem}. 

\begin{lemma}\label{L-Tool02}
    With the notation of Setup \ref{setup} and Remark \ref{R-Observation}, assume that $|E_k \cap E_l| \le 1$ for all $k \ne l.$ For $p \in \Om(\bb)$ and $k \in \Om(\bb, p)$ the following holds.  
   
    \begin{enumerate}[\rm(i)]
        \item $|E_k| = e.$
        \item $E_l \ne E_k$ for all $l \in \Om(\bb, q)$ with $q \in \Om(\bb) \setminus \{p\}.$

        \vspace{2.5mm}
        
        \item  
        $\sum\limits_{j \in E_k \setminus \{p\}} \dfrac{t - r_j}{b_j + 1} \le \dfrac{r_p}{b_p}.$
        
    \end{enumerate}
\end{lemma}

\begin{proof}
    Suppose that $p \in \Om(\bb)$ and $k \in \Om(\bb, p).$ Then, it follows from Remark \ref{R-MonomialOfI}(iii) that $\Lm_k(\bb) = \{p\}$ and $b_p = w_{k, p}.$
    
    \medskip\noindent
    (i) From Remark \ref{R-Observation}(i), we have $a_j = tb_j + r_j$ and $r_j < t$ for all $1 \le j \le n.$ Let $j \in E_k.$ If $j = p$, then 
    \[\dfrac{a_p}{w_{k, p}} = t + \dfrac{r_p}{w_{k, p}} < t + \dfrac{t}{w_{k, p}} \le 2t.\]
    On the other hand, if $j \ne p$, then it follows from Remark \ref{R-MonomialOfI}(i) that $b_j < w_{k, j}$, which implies $a_j = tb_j + r_j < t(b_j + 1) \le t w_{k, j}.$ Thus, it follows from Remark \ref{R-Observation}(i) that
    \[te \le \dfrac{a_p}{w_{k, p}} + \sum_{j \in E_k \setminus \{p\}} \dfrac{a_j}{w_{k, j}} < 2t + (|E_k| - 1)t = (|E_k| + 1)t.\]
    Consequently, $e < |E_k| + 1$, that is $e \le |E_k|.$ Since $e = \max \{|E_i| : 1 \le i \le s\}$, we get $e = |E_k|.$ 
    
    \medskip\noindent
    (ii) Suppose, for contradiction that there exists $q \in \Om(\bb)$ with $q \ne p$, and $l \in \Om(\bb, q)$ such that $E_k = E_l.$ By Remark \ref{R-MonomialOfI}(iii), we have $\Lm_l(\bb) = \{q\}.$ The fact that $q \ne p$ implies that $k \ne l.$ Thus, by given hypothesis that $|E_k \cap E_l| \le 1$, we obtain $|E_k| = |E_l|= 1.$ Therefore, either $\q_k \subset \q_l$ or $\q_l \subset \q_k.$
    This contradicts the fact that $I = \q_1 \cap \cdots \cap \q_s$ is the irredundant irreducible decomposition of $I.$

    \medskip\noindent
    (iii) Using Remark \ref{R-MonomialOfI}(iii), we have $b_p = w_{k, p}.$ Also, it follows from Remark \ref{R-MonomialOfI}(i) that $b_j < w_{k, j}$ for all $j \in E_k$ with $j \ne p.$ Thus, in view of Remark \ref{R-Observation}(i), we have 
    \begin{align*}
        te & \le \sum_{j \in E_k} \dfrac{a_j}{w_{k, j}} = \dfrac{a_p}{w_{k, p}} + \sum_{j \in E_k \setminus \{p\}} \dfrac{a_j}{w_{k, j}} \\  
        & \le \dfrac{tb_p + r_p}{b_p} + \sum_{j \in E_k \setminus \{p\}} \dfrac{tb_j + r_j}{b_j + 1} = t + \dfrac{r_p}{b_p} + \sum_{j \in E_k \setminus \{p\}} \dfrac{t(b_j + 1) + r_j - t}{b_j + 1}.
    \end{align*}
    From part (i), we have $|E_k \setminus \{p\}| = e - 1.$ Thus, we obtain
    \[\sum_{j \in E_k \setminus \{p\}} \dfrac{t(b_j + 1) + r_j - t}{b_j + 1} = t(e - 1) - \sum_{j \in E_k \setminus \{p\}} \dfrac{t - r_j}{b_j + 1}.\]
    Combining these inequalities, we get the desired inequality.
    \[\sum\limits_{j \in E_k \setminus \{p\}} \dfrac{t - r_j}{b_j + 1} \le \dfrac{r_p}{b_p}.\]
\end{proof}

\begin{lemma}\label{L-Tool03}
   With the notation of Set-up \ref{setup} and Remark \ref{R-Observation}, assume that $|E_k \cap E_l| \le 1$ for all $k \ne l.$ If $\Om(\bb) = \emptyset$, then $t\gamma + |\br| \ge t(e - 1).$
\end{lemma}

\begin{proof}
    In view of Remark \ref{R-MonomialOfI}(ii), we have the following two cases.

    \medskip\noindent
    \textit{Case 1.} Assume that $|\Lm_i(\bb)| \ge 2$ for all $1 \le i \le s.$ Let $1\le k \le s$ and $p \in \Lm_k(\bb).$ Then $b_p \ge w_{k, p}$ and $x_p^{w_{k, p}} \in \q_k.$ Consider the monomial 
    \[g = x_p^{w_{k, p}} \prod\limits_{j \notin E_k} x_j^{b_j}.\]
    We first prove that $g \in I.$ It is enough to prove that $g \in \q_i$ for all $1 \le i \le s.$ Clearly, $g \in \q_k.$ Now, let $1 \le i \le s$ with $i \ne k.$ Since $\Lm_i(\bb) \subset E_i$ and $|E_i \cap E_k| \le 1$, there exists $j \in \Lm_i(\bb)$ with $j \notin E_k.$ It follows from Remark \ref{R-MonomialOfI}(ii) that $g \in \q_i.$ This proves that $g \in I.$ 
    
    Now, the fact $b_p \ge w_{k, p}$ implies that $g$ divides $\x^{\bb}.$ As $f \in I$ is a monomial of smallest degree such that $f$ divides $\x^{\bb}$, we obtain $\deg(g) \ge \deg(f).$ This implies that 
    \[|\bb| - \deg(g) \le |\bb| - \deg(f) = \gamma.\]
    From Remark \ref{R-Observation}(i), we have $a_j \ge tb_j$ for all $1 \le j \le n.$ Since $b_p \ge w_{k, p}$, we obtain $g^t$ divides $\x^{\ba}$, and
    \begin{equation}\label{INQ-5}
        \deg\left( \dfrac{\x^{\ba}}{g^t} \right) = t|\bb| + |\br| - t\deg(g) \le t\gamma + |\br|.
    \end{equation}
    One can observe that the monomial 
    \[x_p^{a_p - tw_{k, p}} \prod\limits_{j \in E_k \setminus \{p\}} x_j^{a_j}\]
    divides the monomial $\dfrac{\x^{\ba}}{g^t}.$ Thus, we obtain
    \[\deg\left( \dfrac{\x^{\ba}}{g^t} \right) \ge a_p - tw_{k, p} + \sum_{j \in E_k \setminus \{p\}} a_j \ge \dfrac{a_p - tw_{k, p}}{w_{k, p}} + \sum_{j \in E_k \setminus \{p\}} \dfrac{a_j}{w_{k, j}}.\]
    Thus, it follows from Inequality \ref{INQ-5} and Remark \ref{R-Observation}(i) that
    \[t\gamma + |\br| \ge \dfrac{a_p - tw_{k, p}}{w_{k, p}} + \sum_{j \in E_k \setminus \{p\}} \dfrac{a_j}{w_{k, j}} = \sum_{j \in E_k} \dfrac{a_j}{w_{k, j}} - t \ge t(e - 1).\]

    \medskip\noindent
    \textit{Case 2.} Assume that  $|\Lm_m(\bb)| = 1$ for some $1 \le m \le s.$ Let $\Lm_m(\bb) = \{p\}$, and define 
    \[Z = \big\{ i : 1 \le i \le s \text{ and } \Lm_i(\bb) = \{p\} \big\}.\]
    Choose $k \in Z$ such that 
    \[w_{k, p} = \max \{w_{i, p} : i \in Z\}.\]
    For simplicity, we write $c = b_p - w_{k, p}.$ The fact 
    \[\supp{\x^{\bb}} = \Om(\bb) \sqcup \C(\bb)\]
    and the assumption $\Om(\bb) = \emptyset$ imply $\supp{\x^{\bb}} = \C(\bb).$ Since $p \in \Lm_k(\bb) \subset \supp{\x^{\bb}}$, by using Remark \ref{R-MonomialOfI}(i), we get $b_p - 1 \ge w_{k, p}.$ This implies that $c \ge 1.$ 
    
    Now, consider the partition 
    \[E_k = \{p\} \sqcup A \sqcup B\]
    of $E_k$, where $A = (\supp{\x^{\bb}} \cap E_k) \setminus \{p\}$ and $B = E_k \setminus (A \cup \{p\}).$ Further, let $g$ be the monomial given by
    \[g = \dfrac{\x^\bb}{x_p^{c} \prod\limits_{j \in A} x_j}.\]
    We use Lemma \ref{L-Tool01} to show that $g \in I.$ Clearly, 
    \[\supp{x_p^{c} \prod\limits_{j \in A} x_j} = \{p\} \sqcup A \subset E_k.\] 
    Now, it is enough to prove the following statements.
    \begin{enumerate}[\rm(i)]
        \item $g \in \q_k.$
        
        \item The monomial $h = \dfrac{\x^{\bb}}{x_p^{c}} \in I.$ 

        \item $\dfrac{\x^{\bb}}{x_j} \in I$ for all $j \in A.$ 
    \end{enumerate}
    
    The statement (i) follows from the facts $b_p - c = w_{k, p}$ and that the monomial generator  $x_p^{w_{k, p}} = x_p^{b_p - c}$ of $\q_k$ divides $g.$ 
    
    To prove the statement (ii), let $1 \le i \le s.$ If $i \in Z$, then $w_{i, p} \le w_{k, p} = b_p - c.$ This implies that $x_p^{w_{i, p}}$ divides $h$, and hence $h \in \q_i.$ Now, suppose that $i \notin Z.$ Since $\Lm_i(\bb) \ne \emptyset$, there exists $q \in \Lm_i(\bb)$ such that $q \ne p.$ This implies that $x_q^{b_q}$ divides $h.$ It follows from Remark \ref{R-MonomialOfI}(ii) that $x_q^{b_q} \in \q_i.$ Hence $h \in \q_i.$ This proves statement (ii).
    
    Finally, the statement (iii) follows from the fact that $A \subset \supp{\x^{\bb}} = \C(\bb).$
    
    Thus, by Lemma \ref{L-Tool01}, we get $g \in I.$ Clearly, $g$ divides $\x^{\bb}.$ Since $f \in I$ is a monomial of smallest degree such that $f$ divides $\x^{\bb}$, we obtain $\deg(g) \ge \deg(f).$ Note that 
    \[\deg(g) = |\bb| - (c + |A|).\] 
    This implies that 
    \begin{equation}\label{INQ-6}
        c + |A| = |\bb| - \deg(g) \le |\bb| - \deg(f) = \gamma.
    \end{equation}
    From Remark \ref{R-Observation}(i), we have $a_j = tb_j + r_j$ for all $1 \le j \le n.$ Thus, it follows that 
    \[|\br| \ge \sum_{j \in E_k} r_j = r_p + \sum_{j \in A} r_j + \sum_{j \in B} r_j = (a_p - tb_p) + \sum_{j \in A} (a_j - tb_j) + \sum_{j \in B} (a_j - tb_j).\]
    If $j \in A$, then $j \in E_k$ and $j \ne p.$ Since $\Lm_k(\bb) = \{p\}$, by Remark \ref{R-MonomialOfI}(i), we get $b_j < w_{k, j}.$ Further, if $j \in B$, then $j \notin \supp{\x^{\bb}}$, which means $b_j = 0.$ Also, recall that $c = b_p - w_{k, p}.$ Thus, we obtain
    \begin{align*}
        |\br| & > (a_p - t(c + w_{k, p})) + \sum_{j \in A} (a_j - t w_{k, j}) + \sum_{j \in B} a_j \\
        & \ge - tc + \dfrac{a_p - tw_{k, p}}{w_{k, p}} + \sum_{j \in A} \dfrac{a_j - t w_{k, j}}{w_{k, j}} + \sum_{j \in B} \dfrac{a_j}{w_{k, j}} \\
        & = - t - t(c + |A|) + \sum_{j \in E_k} \dfrac{a_j}{w_{k, j}}.
    \end{align*}
    It follows from Inequality \ref{INQ-6} and Remark \ref{R-Observation}(i) that 
    \[t \gamma + |\br| \ge t(c + |A|) + |\br| > - t + \sum_{j \in E_k} \dfrac{a_j}{w_{k, j}} \ge t(e - 1).\]
\end{proof}

\begin{lemma}\label{L-Tool04}
   With the notation of Set-up \ref{setup} and Remark \ref{R-Observation}, assume that $|E_k \cap E_l| \le 1$ for all $k \ne l.$ If $\Om(\bb) \ne \emptyset$, then $t\gamma + |\br| \ge t(e - 1).$
\end{lemma}

\begin{proof}
    If $\gamma \ge e - 1$, then the result follows immediately. Assume that $\gamma < e - 1.$ Since $\gamma \ge 0$, we have $e > 1.$ To simplify the notation, let $\eta = \max \{b_j : j \in \Om(\bb)\}$, and let $W = \{j \in \Om(\bb) : b_j = \eta\}.$ For each $j \in \Om(\bb)$ and for each $i \in \Om(\bb, j)$, we define the set $A_{j, i}$ as
    \[A_{j, i} = E_i \setminus (\C(\bb) \cup \{j\}).\]
    Then $A_{j, i} \subset E_i \setminus \{j\}.$ First, we prove the following two claims:
    \begin{enumerate}[(i)]
        \item $|A_{j, i}| \ge e - 1 - \gamma > 0.$

        \item $b_{j'} \le \eta$ for all $j' \in A_{j, i}.$ In particular, if $W \cap A_{j, i} = \emptyset$, then $b_{j'} + 1 \le \eta$ for all $j' \in A_{j, i}.$
    \end{enumerate}
    
    We begin by establishing Claim (i). In view of Remark \ref{R-MonomialOfI}(iii) and Lemma \ref{L-Tool02}(i), we get $\Lm_i(\bb) = \{j\}$ and $|E_i| = e.$ By the definition of the sets $\Om(\bb)$ and $\Lm_i(\bb)$, we obtain $j \notin \C(\bb)$ and $j \in E_i.$ Thus, $E_i \cap \C(\bb)$, $\{j\}$, and $A_{j, i}$ are pairwise disjoint subsets of $E_i$, and
    \[E_i = \{j\} \sqcup A_{j, i} \sqcup (E_i \cap \C(\bb)).\]
    One can observe that
    \begin{equation}\label{EQN-Cardinality}
        |A_{j, i}| = e - 1 - |E_i \cap \C(\bb)|.
    \end{equation}
    
    For each $p \in E_i \cap \C(\bb)$, since $\C(\bb) \subset \supp{\x^{\bb}}$, we have $b_p \ge 1.$ Therefore, we can define the monomial
    \[f_j = \dfrac{\x^{\bb}}{\prod\limits_{p \in E_i \cap \C(\bb)} x_p}.\]
    Clearly, the monomial $f_j$ divides $\x^{\bb}.$ We observe the following:
    \vspace{2mm}
    \begin{itemize}
        \item $\supp{\prod\limits_{p \in E_i \cap \C(\bb)} x_p} = E_i \cap \C(\bb) \subset E_i.$

        \vspace{2.5mm}
        
        \item Since $j \notin \C(\bb),$ we get  $x_j^{b_j}$ divides $f_j.$ It follows from Remark \ref{R-MonomialOfI}(ii) and $\Lm_i(\bb) = \{j\}$ that $x_j^{b_j} \in \q_i.$ Consequently, $f_j \in \q_i.$
        
        \vspace{1.5mm}

        \item By using the definition of the set $\C(\bb)$, we obtain $\dfrac{\x^{\bb}}{x_p} \in I$ for all $p \in E_i \cap \C(\bb).$
    \end{itemize}
    Thus, it follows from Lemma \ref{L-Tool01} that $f_j \in I.$ Since $f \in I$ is a monomial of smallest degree such that $f$ divides $\x^{\bb}$, we obtain \[|\bb| - |E_i \cap \C(\bb)| = \deg(f_j) \ge \deg(f).\]
    From Remark \ref{R-Observation}(ii), recall that $\gamma = |\bb| - \deg(f).$ This implies that $\gamma \ge |E_i \cap \C(\bb)|.$ In light of Eqaution \ref{EQN-Cardinality}, we see that
    \[|A_{j, i}| = e - 1 - |E_i \cap \C(\bb)| \ge e - 1 - \gamma.\]
    This proves claim (i).
    
    Now, we prove claim (ii). By the definition of the set $A_{j, i}$, we have $A_{j, i} \cap \C(\bb) = \emptyset.$ By the definition of the sets $\C(\bb)$ and $\Om(\bb)$, we have $\supp{\x^{\bb}} = \C(\bb) \sqcup \Om(\bb).$ Let $j' \in A_{j, i}.$ If $b_{j'} = 0$, then clearly $b_{j'} < \eta.$ On the other hand, if $b_{j'} > 0$, then $j' \in \Om(\bb).$ By the maximality of $\eta$, we get $b_{j'} \le \eta.$ Thus, we obtain $b_{j'} = \eta$ for all $j' \in W \cap A_{j, i}$ and $b_{j'} < \eta$ for all $j' \in A_{j, i} \setminus W.$ This proves claim (ii)

    To establish the proof of the lemma, we choose $p \in W$ such that $r_p = \min \{r_j : j \in W\}$, and let $k \in \Om(\bb, p).$ It follows from Remark \ref{R-MonomialOfI}(iii) that $\Lm_k(\bb) = \{p\}$ and $b_p = w_{k, p}.$ We now consider the following two cases.

    \medskip\noindent
    \textit{Case 1.} Assume that $W \not\subset E_k$, and let $q \in W \setminus E_k.$ Since $\{p\} = \Lm_k(\bb) \subset E_k$, we get $q \ne p.$ Note that the sets $\{p\}$, $A_{p, k}$, and $\{q\}$ are pairwise disjoint. By minimality of $r_p$, we get $r_q \ge r_p \ge \dfrac{r_p}{\eta}.$ Thus, we obtain
    \[|\br| \ge r_p + r_q + \sum_{j \in A_{p, k}} r_j \ge \left( 1 + \dfrac{1}{\eta} \right) r_p + \sum_{j \in A_{p, k}} r_j.\]
    Further, since $p \in W$ and $A_{p, k} \subset E_k \setminus \{p\}$, it follows from claim (ii) and Lemma \ref{L-Tool02}(iii) that
    \[\dfrac{r_p}{\eta} = \dfrac{r_p}{b_p} \ge \sum_{j \in E_k \setminus \{p\}} \dfrac{t - r_j}{b_j + 1} \ge \sum_{j \in A_{p, k}}\dfrac{t - r_j}{b_j + 1} \ge \dfrac{1}{\eta + 1} \left( t |A_{p, k}| - \sum_{j \in A_{p, k}} r_j \right).\]
    Applying claim (i), we obtain
    \[|\br| \ge \left( 1 + \dfrac{1}{\eta} \right) r_p + \sum_{j \in A_{p, k}} r_j \ge t |A_{p, k}| \ge t(e - 1 - \gamma).\]

    \medskip\noindent
    \textit{Case 2.} Assume that $W \subset E_k.$ We first claim that there exists $q \in W$ such that $W \cap E_l = \{q\}$ for all $l \in \Om(\bb, q).$ 
    
    We first consider the case  $W = \{p\}.$ For each $l \in \Om(\bb, p)$, by Remark \ref{R-MonomialOfI}(iii), we have $\Lm_l(\bb) = \{p\}.$ Since $\Lm_l(\bb) \subset E_l$, we obtain $W \cap E_l = \{p\}$ for all $l \in \Om(\bb, p).$ Thus, if $W = \{p\}$, then we can take $q = p.$ Now, suppose that $\{p\} \subsetneq W$, and let $q \in W$ with $q \ne p.$ If $l \in \Om(\bb, q)$, then by Remark \ref{R-MonomialOfI}(iii), we have $\Lm_l(\bb) = \{q\}.$ Since $\Lm_l(\bb) \subset E_l$ and $q \in W$, we get $\{q\} \subset W \cap E_l.$ Also, from Lemma \ref{L-Tool02}(ii), we have $E_l \ne E_k.$ Thus, by the given hypothesis 
    \[|W \cap E_l| \le |E_k \cap E_l| \le 1.\]
    This shows that $W \cap E_l = \{q\}.$

    Now, let $q \in W$ be as above, and let $l \in \Om(\bb, q)$ be fixed. Then $b_q = \eta$ and $W \cap E_l = \{q\}.$ Since $A_{q, l} \subset E_l \setminus \{q\}$, it follows that $W \cap A_{q, l} = \emptyset.$ By claim (ii) and Lemma \ref{L-Tool02}(iii), we have
    \[\dfrac{r_q}{\eta} = \dfrac{r_q}{b_q} \ge \sum_{j \in E_l \setminus \{q\}} \dfrac{t - r_j}{b_j + 1} \ge \sum_{j \in A_{q, l}}\dfrac{t - r_j}{b_j + 1} \ge \dfrac{1}{\eta} \left( t |A_{q, l}| - \sum_{j \in A_{q, l}} r_j \right).\]
    The fact $A_{q, l} \subset E_l \setminus \{q\}$ implies that $q \notin A_{q, l}.$ Applying claim (i), we obtain
    \[|\br| \ge r_q + \sum_{j \in A_{q, l}} r_j \ge t |A_{q, l}| \ge t(e - 1 - \gamma).\]
    This completes the proof of the lemma.
\end{proof}

\begin{thm}\label{Th-MainTheorem}
    With the notation of Setup \ref{setup} and Remark \ref{R-Observation}, assume that $|E_k \cap E_l| \le 1$ for all $k \ne l.$ Then, for every integer $t \ge 1$, we have the containment
    \[I^{(te)} \subset \mi^{t(e - 1)} I^{[t]}.\]
\end{thm}

\begin{proof}
    In view of Remark \ref{R-Observation}(ii), it is enough to prove that $t \gamma + |\br| \ge t(e - 1).$ This assertion follows from Lemmas \ref{L-Tool03} and \ref{L-Tool04}.
\end{proof}

The following corollary is a direct consequence of Theorem \ref{Th-MainTheorem}.

\begin{corollary}\label{C-LowerBound}
    With the notation of Setup \ref{setup} and Remark \ref{R-Observation}, assume that $|E_k \cap E_l| \le 1$ for all $k \ne l.$ Then, for every integer $t \ge 1$, we have the containment
    \[I^{(te)} \subset \mi^{t(e - 1)} I^t.\]
    Moreover, we have
    \[\ww(I) \ge \dfrac{\alpha(I) + e - 1}{e}.\]
\end{corollary}

\begin{thm}\label{Th-RevisedBound}
    Let $I \subset R$ be a monomial ideal of big-height $e$, and let
    \[w_0 = \min\{\deg_{x_j}(f) : f \in \mathcal{G}(I) \text{ and } j \in \supp f\}.\] Then the containment $I^{(t(e + 1))} \subset \mi^{[tw_0]} I^{[t]}$ holds for all $t \ge 1.$
\end{thm}

\begin{proof}
    We follow the notation of Set-up \ref{setup}. Note that 
    \[w_0 = \min\{w_{i, j} : 1 \le i \le s \text{ and } j \in E_i\}.\]
    Let $\x^{\ba} = x_1^{a_1} \cdots x_n^{a_n} \in I^{(t(e+1))}.$ Then it follows from Lemma \ref{L-SymbolicPowerInequality} that 
    \begin{equation}\label{INQ-0}
        \sum_{j \in E_i} \dfrac{a_j}{w_{i, j}} \ge t(e+1).
    \end{equation}
    Note that $e = \max\{ |E_i| : 1 \le i \le s\}.$ Thus, for each $1 \le i \le s$, there exists $j_i \in E_i$ such that 
    \[\dfrac{a_{j_i}} {w_{i, {j_i}}} > \dfrac{t(e + 1)}{|E_i|} \ge \dfrac{t(e + 1)}{e} > t.\]
    This implies that $a_{j_i} > t w_{i, {j_i}} \ge t w_0$ for all $1 \le i \le s.$ Now, consider the monomial $f = \dfrac{\x^{\ba}}{x_{j_1}^{t w_0}}.$ Since $x_{j_1}^{t w_0} \in \mi^{[tw_0]}$, it suffices to prove that $f \in I^{[t]}$, or equivalently, $f \in \q_i^{[t]}$ for all $1 \le i \le s.$ Let $1 \le i \le s.$ We consider the following two cases.
    
    \medskip\noindent
    \textit{Case 1.} When $j_1 \in E_i.$ If $\dfrac{a_j}{w_{i, j}} \ge t$ for some $j \in E_i$ with $j \ne j_1$, then we obtain
    \[\deg_{x_j}(f) = \deg_{x_j}(\x^{\ba}) = a_j \ge tw_{i, j}.\]
    This implies that $f \in \q_i^{[t]}.$ On the other hand, if $\dfrac{a_j}{w_{i, j}} < t$ for all $j \in E_i$ with $j \ne j_1$, then we have
    \[\sum_{j \in E_i \setminus \{j_1\}} \dfrac{a_j}{w_{i, j}} < t(|E_i| - 1) \le t(e - 1).\]
    In view of Inequality \ref{INQ-0}, we obtain
    \[\dfrac{a_{j_1}}{w_{i, j_1}} = \sum_{j \in E_i} \dfrac{a_j}{w_{i, j}} - \sum_{j \in E_i \setminus \{j_1\}} \dfrac{a_j}{w_{i, j}} \ge t(e + 1) - t(e - 1) = 2t.\]
    This implies that $a_{j_1} \ge 2t w_{i, j_1}$, which means $\deg_{x_{j_1}}(f) =  a_{j_1} - t w_0 \ge t w_0.$ Hence $f \in \q_i^{[t]}.$

    \medskip\noindent 
    \textit{Case 2.} When $j_1 \notin E_i.$ In this case, we have $\deg_{x_{j_i}}(f) = a_{j_i}.$ Since $\dfrac{a_{j_i}}{w_{i, j_i}} > t$, it follows that $\deg_{x_{j_i}}(f) = a_{j_i} > t w_{i, j_i}.$ This implies that $f \in \q_i^{[t]}.$
\end{proof}

\begin{corollary}\label{C-RevisedBound}
    With the notation of Set-up \ref{setup}, let $w_0 = \min \left\{w_{i, j} : 1 \le i \le s, \text{ and } j \in E_i \right\}.$ If $e$ is the big-height of $I$, then
    \[\ww(I) \ge \dfrac{\alpha(I) + w_0}{e + 1}.\]
    In particular, if $w_0 \ge e - 1$, then 
    \[\ww(I) \ge \dfrac{\alpha(I) + e - 1}{e + 1}.\]
\end{corollary}

The aforementioned result naturally lead to the following question. 

\begin{question}
    Let $I \subset R$ be a monomial ideal of big-height $e.$ Is the lower bound \[\ww(I) \ge \dfrac{\alpha(I) + e - 1}{e + 1}\] always satisfied? 
\end{question}

\section*{Acknowledgments}

The authors acknowledge the use of the computer algebra system Macaulay2 \cite{M2}  for testing their computations. The first author is supported by the Indian Institute of Technology Jammu, under the institute fellowship.

\subsection*{Conflict of interest} The authors declare that they have no known competing financial interests or personal relationships that could have appeared to influence the work reported in this paper.

\section* {AI Disclosure}

OpenAI’s ChatGPT free version was used as an auxiliary tool to identify and correct typographical and grammatical errors throughout the manuscript. The AI was consulted exclusively in connection with certain mathematical arguments involved in the proof of Theorem \ref{Th-MainTheorem}. Initially, AI failed to produce any mathematically valid arguments, and several subsequent responses were inconclusive. The authors then provided additional mathematical directions and guidance, following which AI generated conclusive arguments that assisted the authors in formulating a rigorous proof of Theorem \ref{Th-MainTheorem}. The authors critically examined, verified, and refined these arguments before incorporating them into the final proof. AI assistance was limited exclusively to the proof of Theorem \ref{Th-MainTheorem}. All other mathematical results and proofs presented in this manuscript, including those in Section \ref{S-CounterExample} as well as Theorem \ref{Th-RevisedBound} and Corollary \ref{C-RevisedBound}, are developed independently, without the assistance of AI. The authors take full responsibility for the mathematical statements, the validity and correctness of all proofs, and the final content of the manuscript.

\bibliographystyle{amsplain}
\bibliography{references}

\end{document}